\documentclass[12pt]{amsart}
\usepackage[margin=1.5in]{geometry}
\usepackage{amsthm,amsmath,amssymb}
\theoremstyle{plain}
\newtheorem{theorem}{Theorem}[section]

\newtheorem{question}[theorem]{Question}
\theoremstyle{definition}

\numberwithin{equation}{section}
\newcommand*{\abs}[1]{\lvert#1\rvert}
\newcommand*{\set}[1]{\{#1\}}
\newcommand*{\defset}[2]{\{\,#1\,\mid\,#2\,\}}
\newcommand{\smalldet}[4]{\begin{vmatrix}#1 & #2\\#3 &#4\end{vmatrix}}
\newcommand*{\QQ}{\mathbb Q}
\newcommand*{\ZZ}{\mathbb Z}

\newcommand*{\RR}{\mathbb R}
\newcommand*{\FF}{\mathbb F}
\DeclareMathOperator{\SL}{SL}
\usepackage{xcolor} 
\usepackage{url}
\usepackage[T1]{fontenc}
\usepackage{lmodern}
\usepackage{listings}
\definecolor{kw}{RGB}{0,0,150}
\definecolor{str}{RGB}{150,0,0}
\definecolor{cmt}{RGB}{60,120,60}

\lstdefinelanguage{Sage}{
  language=Python,
  morekeywords={matrix,GF,identity,diagonal}
}
\usepackage[
  backend=biber,
  style=numeric,
  doi=true,
  url=false,
  giveninits=true,
  isbn=false
]{biblatex}
\usepackage[colorlinks,linkcolor=blue,citecolor=purple]{hyperref}
\title{Explicitly combing hedgehogs over fields of Stufe $4$}
\author{Peter M\"uller} \address{Institute of Mathematics, University
  of W\"urzburg} \email{peter.mueller@uni-wuerzburg.de}
\begin{document}
\begin{abstract}
  Let
\[
  K[x,y,z]=K[X,Y,Z]/(X^2+Y^2+Z^2-1)
\]
be the coordinate ring of the algebraic unit sphere over a field
$K$. Umberto Zannier showed that there exists a matrix in
$\SL_3(K[x,y,z])$ with first row $(x,y,z)$ for $K=\QQ_p$, the field of
$p$-adic numbers for an odd prime $p$, or more generally, if $-1$ is a
sum of two squares in $K$.
The case $K=\QQ_2$ remained open and was subsequently posed and
discussed by Zannier with numerous researchers, thereby bringing the
problem to broader attention.

In \cite{AnaLev:combing}, Alexey Ananyevskiy and Marc Levine showed
that such a matrix exists if and only if $K$ has Stufe at most $4$,
equivalently, if there exist $a,b,c,d\in K$ such that
$a^2+b^2+c^2+d^2=-1$. Since $\QQ_2$ has Stufe $4$, this settled
Zannier's problem.

Their proof is purely existential and does not provide an explicit
matrix. In this note, we construct an explicit example in terms of
$a,b,c,d$ and describe the computational techniques used to find it.
\end{abstract}
\maketitle
\section{Introduction}
\subsection{The unimodular matrix completion problem of Zannier}
Let $X$, $Y$, $Z$ be variables over a field $K$ and $x, y, z$ be the
images of $X, Y, Z$ in the coordinate ring $K[X,Y,Z]/(X^2+Y^2+Z^2-1)$
of the algebraic unit sphere over $K$. Umberto Zannier gave an
unpublished positive answer to the following question for all $p$-adic
fields $K=\QQ_p$ for an odd prime $p$.
\begin{question}\label{que}
  Does there exist a matrix in $\SL_3(K[x,y,z])$ with first row
  $(x,y,z)$?
\end{question}
The \emph{Stufe} of a field $K$ is defined to be $\infty$ if $-1$ is
not a sum of squares in $K$, and otherwise the smallest number of
squares whose sum equals $-1$. A classical theorem of Pfister states
that every finite Stufe is a power of $2$.

The fields $\QQ_p$ for odd $p$ have Stufe $2$. This is what Zannier
used. For if $a^2+b^2=-1$ with $a,b\in K$, then 
\[
  \left(\begin{array}{ccc}
    x & y & z \\
    1 & a & b \\
    0 & b x + z & -a x - y
  \end{array}\right)\in\SL_3(K[x, y, z])
\]
is a positive example to Question \ref{que}.

The field $\QQ_2$, however, has Stufe $4$, so this simple argument
does not apply. Zannier subsequently posed the question for $K=\QQ_2$
to numerous mathematicians working in number theory, algebraic
geometry, and algebraic $K$-theory, helping bring it to broader
attention. He was moreover interested, should such matrices exist, in
explicit constructions rather than merely existential arguments.

Eventually, in \cite{AnaLev:combing}, Alexey Ananyevskiy and Marc
Levine settled Zannier's question. In fact they proved that Question
\ref{que} has a positive answer if and only if $K$ has Stufe at most
$4$.

Their abstract proof, which uses quite a bit of advanced methods and
results, is purely existential and does not produce an explicit
matrix.

We complement their theoretical existence result by an explicit
example:
\begin{theorem}\label{theo}
  Let $K$ be a field, $a^2+b^2+c^2+d^2=-1$ for $a, b, c, d\in K$, and
  $K[x, y, z]=K[X, Y, Z]/(X^2+Y^2+Z^2-1)$ as above.  For
\begin{align*}
m_1 &= (-2 a b c + a^{2} - b^{2}) y + (a^{2} c - b^{2} c + 2 a b) z\\
m_2 &= (a b c - a^{2} + c^{2} + 1) x + (b c^{2} + a c + 2 b) y + (-a c^{2} + b c - 2 a) z - a c d\\
m_3 &= (-a^{2} c + c^{3} - a b + c) x + (2 a c^{2} + b c + a) y + (2 b c^{2} - a c + b) z + a d\\
m_4 &= (-a^{2} b d - b^{3} d) z + 2 a^{2} b c - a^{3} + a b^{2}\\
m_5 &= (a b^{2} d + b c d - a d) x + (b^{2} c d - 2 a b d) y + (-a b c
      d + a^{2} d + b^{2} d) z\\
  &\phantom{==}- a^{3} b + a b^{3} + a b c^{2} - a^{2} c\\
m_6 &= (a^{2} b d + b c^{2} d - a c d) x + (-a^{2} d) y + (2 b^{2} c d - a b d) z + 2 a^{2} b^{2} - a b c + a^{2}
\end{align*}
set\[
  M = \left(\begin{matrix}x & y & z\\m_1 & m_2 & m_3\\m_4 & m_5 & m_6
  \end{matrix}\right).\]
Then \[
  \det(M)=2a^2(a^2 + d^2)(2abc - a^2 + b^2).
  \]
\end{theorem}
Note that this theorem provides an explicit answer to Question
\ref{que} for all fields $K$ of Stufe $\le4$. We only need to make
sure that $\det(M)$ does not vanish. We consider three cases:
\begin{itemize}
\item[(a)] If $K$ has Stufe $1$, then set $a=\sqrt{-1}$, $b=c=d=0$.
  \item[(b)] If $K$ has Stufe $2$, then there are $a,d\in K$ with
    $a^2+d^2=-1$ and $a\ne0$. Set $b=c=0$.
  \item[(c)] If $K$ has Stufe $4$, then pick $a,b,c,d\in K$ with
    $a^2+b^2+c^2+d^2=-1$. Then $a,b,c,d$ are nonzero, and so is
    $a^2+d^2$, for otherwise the Stufe of $K$ would be $<4$. Thus if
    $\det(M)=0$, then $2abc=a^2-b^2$. Now switch $a$ with $b$. If the
    corresponding determinant would vanish again, then $2abc=b^2-a^2$,
    so $abc=0$, a contradiction.
  \end{itemize}
  Zannier had also asked whether Question \ref{que} has a positive
  answer for $K=\ZZ_2$, the $2$-adic integers. While
  \cite{AnaLev:combing} handles $\QQ_2$, it does not help for this
  sharpened question. Nevertheless, we get an example over $\ZZ_2$
  from setting
  $(a, b, c, d) = (\frac{\omega}{7}, 1, \frac{3\omega}{7},
  \frac{2\omega}{7})$ in Theorem \ref{theo}, where $\omega=\sqrt{-7}$:
\begin{theorem}\label{theo7}
  For $\omega=\sqrt{-7}$ consider the matrix
  \[
    M=\left(\begin{array}{ccc}
      x & y & z \\
      \omega y - 5 z & 2 \omega x - \omega y + 8 z + 3 & -5 x + 5 y + 4 \omega z + \omega \\
      -6 z + 1 & 3 \omega x + \omega y + 9 z + 1 & -7 x + y + 5 \omega
                                                   z
\end{array}\right)
\]
with entries in
$\ZZ[\omega][x, y, z]=\ZZ[\omega][X, Y, Z]/(X^2+Y^2+Z^2-1)$. Then
$\det(M)=5$.
\end{theorem}
Note that $\ZZ[\frac{1}{5},\omega]\subset\ZZ_2$, so upon dividing the
second row of $M$ by $5$, the matrix is still in
$\ZZ_2[x,y,z]^{3\times 3}$. A similar example was communicated to the authors of
\cite{AnaLev:combing} by the author of the present note, see \cite[Remark
1.10]{AnaLev:combing}.
\subsection{Nowhere vanishing vector fields on spheres.}
A matrix $M$ as in Question \ref{que} gives a
nowhere vanishing vector field on the algebraic sphere in the
following sense. Let $\bar K$ be an algebraic closure of $K$ and
assume that the characteristic of $K$ is not $2$. Let
\[
S(\bar K)=\{(x_0,y_0,z_0)\in \bar K^3 : x_0^2+y_0^2+z_0^2=1\}
\]
be the $\bar K$-points of the algebraic $2$-sphere over $K$, and let
$w(x, y, z)$ be the second row of $M\in\SL_3(K[x, y, z])$. As
$x^2+y^2+z^2=1$, we see that $w(x, y, z)$ defines a well-defined
polynomial map on $S(\bar K)$. As $\det(M)$ is a nonzero constant in
$K$, we obtain that the vectors $(x_0, y_0, z_0)$ and
$w(x_0, y_0, z_0)$ are linearly independent for all
$(x_0, y_0, z_0)\in S(\bar K)$.

With $\cdot$ the dot product of vectors set
\[
  v(x, y, z)=w(x, y, z)-((x, y, z)\cdot w(x, y, z))(x,y,z).
\]
Then $v(x_0, y_0, z_0)\cdot(x_0, y_0, z_0)=0$ and $v(x_0, y_0,
z_0)\ne0$ for all $(x_0, y_0, z_0)\in S(\bar K)$, so $v(x,y,z)$
defines a nowhere vanishing vector field on the sphere $S(\bar K)$.

This explains the title of this note.

Conversely, from a nowhere vanishing vector field on $S(\bar K)$ one
can construct a matrix as in Question \ref{que}, see Appendix
\ref{app}.
\subsection{Nowhere vanishing vector fields on $S(\bar K)$ and on
  $S(K)$.}
By the Hairy Ball Theorem from real topology, there is no nowhere
vanishing vector field on $S(\RR)$. So in particular, Question
\ref{que} has a negative answer for $K=\RR$.

In view of the real case, an alternative $p$-adic analog of the
``combing the hedgehog'' question would only require that a vector
field $v(x, y, z)$ does not vanish on $S(K)$. In this case, there is
an easy example even over $\QQ_2$: Set $v(x, y, z)=(y+z, -x,
-x)$. Then $(x, y, z)\cdot v(x, y, z)=0$, and $v(x, y, z)=0$ would
force that $2$ is a square in $\QQ_2$ (recall that $x^2+y^2+z^2=1$),
which is not the case.
\subsection{Acknowledgement} I am grateful to Umberto Zannier for
formulating and widely disseminating his question that led to this
note, for several illuminating discussions about it, and for helpful
comments on an earlier draft.

I thank Alexey Ananyevskiy and Marc Levine for helpful discussions and
for answering questions about their paper.

I also thank Joachim König for telling me about Zannier's problem,
correctly anticipating that I would find it attractive.
\section{Verification of Theorems \ref{theo} and \ref{theo7}}\label{proof}
The following \texttt{SageMath} \cite{sagemath} code verifies Theorems
\ref{theo} and \ref{theo7}:
\begin{lstlisting}
#### Verifying Theorem 1.2 ####

R.<a, b, c, d, X, Y, Z> = QQ[]
           
m1 = (-2*a*b*c + a^2 - b^2)*Y + (a^2*c - b^2*c + 2*a*b)*Z
m2 = ((a*b*c - a^2 + c^2 + 1)*X + (b*c^2 + a*c + 2*b)*Y
      + (-a*c^2 + b*c - 2*a)*Z - a*c*d)
m3 = ((-a^2*c + c^3 - a*b + c)*X + (2*a*c^2 + b*c + a)*Y
      + (2*b*c^2 - a*c + b)*Z + a*d)
m4 = (-a^2*b*d - b^3*d)*Z + 2*a^2*b*c - a^3 + a*b^2
m5 = ((a*b^2*d + b*c*d - a*d)*X + (b^2*c*d - 2*a*b*d)*Y
      + (-a*b*c*d + a^2*d + b^2*d)*Z - a^3*b + a*b^3 + a*b*c^2 - a^2*c)
m6 = ((a^2*b*d + b*c^2*d - a*c*d)*X + (-a^2*d)*Y
      + (2*b^2*c*d - a*b*d)*Z + 2*a^2*b^2 - a*b*c + a^2)

M = matrix([[X, Y, Z], [m1, m2, m3], [m4, m5, m6]])
detM = 2 * a^2 * (a^2 + d^2) * (a^2 - b^2 - 2*a*b*c)

q1 = X^2 + Y^2 + Z^2 - 1
q2 = a^2 + b^2 + c^2 + d^2 + 1

assert q2.divides((M.det() - detM)%q1)

#### Verifying Theorem 1.3 ####

K.<w> = QuadraticField(-7)
R.<X, Y, Z> = K[]

M = matrix([
    [X, Y, Z],
    [w*Y - 5*Z, (2*w)*X + (-w)*Y + 8*Z + 3, -5*X + 5*Y + (4*w)*Z + w],
    [-6*Z + 1, (3*w)*X + w*Y + 9*Z + 1, -7*X + Y + (5*w)*Z]])

assert M.det()%(X^2 + Y^2 + Z^2 - 1) == 5
\end{lstlisting}
The code is also available at \cite{verify_combing}, where it can be
downloaded or run online.
\section{The search for the example in Theorem \ref{theo}}\label{comp}
Finding the example in Theorem \ref{theo} required quite some
effort. As we believe that the techniques and ideas can be useful for
similar questions, we describe the main ideas and tools.

For the computations described below we used \texttt{SageMath}
\cite{sagemath}, and for difficult Groebner basis computations,
including computing elimination ideals and rational univariate
representations of solution sets of $0$-dimensional varieties, we used
the excellent software \texttt{msolve} \cite{msolve}.

Let $K$ be a field of Stufe $4$. So in particular, the characteristic
of $K$ is $0$. We are looking for a matrix \[ M=%
  \left(\begin{matrix}X & Y & Z\\m_1 & m_2 & m_3\\m_4 & m_5 & m_6
    \end{matrix}\right)\in K[X, Y, Z]^{3\times 3}\]
such that $\det(M)\equiv1\pmod{X^2+Y^2+Z^2-1}$.

One first confirms that in contrast to the case where the field $K$
has Stufe $\le2$, each row of $M$ needs to contain a non-constant
polynomial: Suppose otherwise. Then without loss of generality,
$m_1=1$, $m_2=a\in K$, $m_3=b\in K$. Pick $\gamma\in\bar K$ with
$\gamma^2=a^2+b^2+1$. Then $\gamma\ne0$, because $K$ has Stufe
$4$. Set
\[ (x, y, z) = (\frac{1}{\gamma}, \frac{a}{\gamma},
  \frac{b}{\gamma}).
\]
Then $x^2+y^2+z^2=1$, and the first two rows of the matrix
$M(x, y, z)$ are linearly dependent, a contradiction.

We resorted to the optimistic assumption that all entries of $M$ have
degree at most $1$. As $\det(M)$ has degree $\le3$ in this case, the
condition \[\det(M)\equiv1\pmod{X^2+Y^2+Z^2-1}\] translates to
\[
  \det(M)=1+(b_0X+b_1Y+b_2Z+b_3)(X^2+Y^2+Z^2-1),
\]
with unknowns $b_0, b_1, b_2, b_3$ over $K$.

Set
{
  \scriptsize
\[
M=\left(\begin{array}{ccc}
X & Y & Z \\
a_{0} X + a_{1} Y + a_{2} Z + a_{3} & a_{4} X + a_{5} Y + a_{6} Z + a_{7} & a_{8} X + a_{9} Y + a_{10} Z + a_{11} \\
a_{12} X + a_{13} Y + a_{14} Z + a_{15} & a_{16} X + a_{17} Y + a_{18}
                                          Z + a_{19} & a_{20} X +
                                                       a_{21} Y +
                                                       a_{22} Z +
                                                       a_{23}
\end{array}\right)
\]
}for unknowns $a_0, a_1,\dots,a_{23}$. Comparing coefficients
instantly expresses $b_0$, $b_1$, $b_2$, $b_3$ in terms of the
$a_i$'s.

This gives a set $Q\subset\ZZ[a_0,a_1,\dots,a_{23}]$ of $16$ quadratic
polynomials in the $a_i$'s such that the matrices $M$ as above
correspond to the simultaneous solutions of $Q$ in $K^{24}$. More
generally, if $R$ is a ring, we write
\[
  V(R)=\defset{(\bar a_0,\dots,\bar a_{23})\in R^{24}}%
  {f(\bar a_0,\dots,\bar a_{23})=0\text{ for all }f\in Q}
\]
for the set of $R$-solutions of $Q$. In the following $R$ will be
either $K$, or the $2$-adic integers $\ZZ_2$, or the field $\FF_2$
with two elements.

Pick $(\bar a_0,\dots,\bar a_{23})\in V(R)$. As $Q$ contains the polynomial
\[
  a_{11} a_{13} + a_{9} a_{15} - a_{3} a_{21} - a_{1} a_{23} - 1,
\]
we see that at least one of $\bar a_{11}$, $\bar a_{15}$,
$\bar a_{3}$, or $\bar a_{23}$ is nonzero. By switching (and
rescaling) rows and columns and permuting $X, Y, Z$, we may and do
assume that $\bar a_{15}$ is nonzero.

If $R=K$, then dividing the third row of $M$ by $\bar a_{15}$ and
multiplying the second row by $\bar a_{15}$ allows us to assume $\bar
a_{15}=1$. This may also be assumed if $R=\ZZ_2$ with essentially the
same argument, because at least one of $\bar a_{11}$, $\bar a_{15}$,
$\bar a_{3}$, $\bar a_{23}$ maps to $1$ in $\FF_2$, so is a unit of
$\ZZ_2$.

Thus we may assume $a_{15}=1$.

With $a_{15}=1$, we may kill $a_{3}$ in the second row, while
preserving all assumptions.

There is one further simplification, for by adding multiples of the
first row to the other rows we may assume $a_0=a_{12}=0$.

Thus so far we may and do assume $a_0=a_3=a_{12}=0$ and $a_{15}=1$.

Now we look at $2$-adic solutions $V(\ZZ_2)$. Note that each such
solution maps to an element in $V(\FF_2)$.

Computing $V(\FF_2)$ isn't hard. One option is to add the polynomials
$a_i^2-a_i$ to $Q$, compute a Groebner basis in lexicographical term
order, and solve the system over $\FF_2$. Here the easier option is to
just run through the $2^{20}$ candidate tuples, which takes a few
seconds only when naively implemented in Python. As a result,
$\abs{V(\FF_2)}=80$.

I wrote a Sage program which does a depth first backtracking search
for lifts of $\FF_p$-solutions to integers solutions modulo $p^k$.

Among these $80$ $\FF_2$-solutions, $76$ don't even lift to integer
solutions modulo $4$.

However, the remaining four $\FF_2$-solutions lift to integer
solutions modulo $2^{50}$, which makes it very likely that they lift
to $\ZZ_2$-solutions. In these four cases $\bar a_1$ is odd. As we
already fixed $a_{15}=1$ in the other row, we cannot assume
$a_1=1$. However, as $\bar a_1$ is a unit in $\ZZ_2$, we can kill
$a_{13}$ in the other row.

In the case $R=K$, there is no justification to assume $\bar a_1\ne0$,
except for the expectation that generically $\bar a_1$ should not
vanish. That is, forcing $a_1=0$ would yield only a small subset of
$V(K)$.

So eventually we assume $a_0=a_3=a_{12}=a_{13}=0$ and $a_{15}=1$.

At this point, the ``without loss of generalities'' are exhausted. So
we have $16$ polynomial conditions for $24-5=19$ unknowns.

The ideal generated by $Q$ is too complex to compute a Groebner
basis. So in particular, we cannot directly compute its
dimension. However, if $V(\bar K)$ is not empty, then the dimension is
at least $3$.

Initially, we expected the dimension to be exactly $3$. So
specializing three of the unknowns to rationals should in general give
a $0$-dimensional system. And indeed, this happened for specializing
certain triples of variables, and the solutions usually lived in
$\QQ_2$.

However, for most triples of the variables, specializing them yielded
a system of positive dimension. Only by specializing one more
variable, we usually got a $0$-dimensional system. But in all tested
cases none of the solutions lived in $\QQ_2$.

This observation suggested that the variety $V(\bar K)$ has at least
two components, where one component has dimension $3$ and solutions in
$\QQ_2$, and the other component has dimension $4$ and no
$\QQ_2$-solutions.

Of course, as the unspecialized system $Q$ is too large to compute a
Groebner basis, it is even more impossible to compute a primary
decomposition.

However, after specializing $4$ of the variables in the cases where
specializing only $3$ of them results in a system of positive
dimension, \texttt{msolve} could compute Groebner bases. And it turned
out that in these cases $a_{19}^2+a_{23}^2+1$ is in the ideal
generated by $Q$.

But $a_{19}^2+a_{23}^2+1=0$ has no solution in a field with Stufe $4$.

Thus, we extended our list of variables by the slack variable $h$ and
$Q$ by the polynomial $h\cdot(a_{19}^2+a_{23}^2+1)-1$.

And indeed, randomly specializing three variables usually resulted in
$0$-dimensional varieties with points in $\QQ_2$.

In particular, rational specializations of $a_{17}$, $a_{21}$,
$a_{22}$ always resulted in solutions over quadratic extensions of
$\QQ$.

Thus, for $0\le i\le 23$ with
$i\not\in\set{0, 3, 12, 13, 15, 17, 21, 22}$, there should be a
polynomial in $a_i,a_{17},a_{21},a_{22}$ of degree at most $2$ in
$a_i$, and in cases where the degree is $2$, the squarefree part of
the discriminant with respect to $a_i$ should be independent of $i$.

Usually, one can obtain such relations by computing elimination
ideals. However, various Groebner packages, including \texttt{msolve},
terminated after using 128GB of memory!

So we tried a different strategy. First note that if we know all
$\bar a_i$, $i\ge12$, of a solution, then determining $\bar a_i$,
$i<12$ amounts to solving a system of linear equations.

Thus we were trying to compute low degree polynomial relations between
$a_i,a_{17},a_{21},a_{22}$, $i\ge12$, by computing low degree
polynomials in $a_i,a_{17},a_{21},a_{22}$ which vanish in thousands of
solutions over quadratic fields. Note that this amounts in solving
linear equations. However, except for $i=18$, where it happened that
$a_{18}$ depends linearly on $a_{17},a_{21},a_{22}$, there was no
polynomial relation of degree $\le10$.

Next we tried the following: We computed a basis of the $\QQ$-vector
space of degree $\le4$ polynomials in the nine variables $a_i$,
$14\le i\le 23$, $i\ne15$ which vanish in thousands of solutions. The
vector space of these polynomials has dimension $306$, and
\texttt{msolve} was able to compute elimination ideals of the ideal
generated by the $306$ basis vectors. In particular, the $a_i$'s we
are looking for, i.e.\ for $i\in\set{14, 16, 19, 20, 23}$, depend only
quadratically on $a_{17},a_{21},a_{22}$. (If we had used a basis of
degree $\le3$ relations, then the dependency of the $a_i$'s would have
been of degree $4$.)

The lowest degree relation is between the variables
$a_{14}, a_{17},a_{21},a_{22}$ with $f(a_{14},a_{17},a_{21},a_{22})=0$
for a degree $10$ polynomial $f$ over $\QQ$. The individual degrees
with respect to $a_{14},a_{17},a_{21},a_{22}$ are $2, 8, 6, 8$,
respectively.

While $f$ is quite a complicated polynomial, it turned out that the
substitution
\begin{align*}
  a_{14} &= r(4s^2 + t^2 - 4s(2s-u) + (2s-u)^2)/t\\
  a_{17} &= s\\
  a_{21} &= t/3\\
  a_{22} &= 2s-u
\end{align*}
for new variables $r, s, t, u$ considerably simplifies $f$. We get $f(a_{14},a_{17},a_{21},a_{22})=(t^2 + u^2)^2F(r, s, t, u)$,
where
\begin{multline*}
  F(r, s, t, u)=36 r^{2} s^{2} t^{4} - 36 r^{2} s t^{4} u + 9 r^{2}
  t^{4} u^{2} - 12 r^{2} s t^{2} u^{3} + 6 r^{2} t^{2} u^{4} + r^{2}
  u^{6} + 9 r^{2} t^{2} u^{2}\\ + 9 r^{2} u^{4} + 9 s^{2} t^{2} - 12 s
  t^{2} u + 4 t^{2} u^{2} + u^{4}.
\end{multline*}
Note that $(t^2 + u^2)^2$ cannot not vanish, because $t\ne0$ and the
field has Stufe $4$. Thus for a solution, $F$ has to vanish. However,
as $V(\RR)=\set{}$ (because the hedgehog cannot be combed over the
reals), $F=0$ does not have a real solution. Together with
$F(0,0,0,1)=1>0$, we see that $F$ is positive on $\RR^4$. By Artin'
Theorem, $F$ is a sum of squares of rational functions. Here, however,
the situation is much better. $F$ is actually a sum of five squared
polynomials. We have $F=q_0^2+q_1^2+q_2^2+q_3^2+q_4^4$, where
\begin{align*}
  q_0 &= u^{2}\\
  q_1 &= -3 s t + 2 t u\\
  q_2 &= 3 r t u\\
  q_3 &= 3 r u^{2}\\
  q_4 &= 6 r s t^{2} - 3 r t^{2} u - r u^{3}.
\end{align*}
Thus the condition $F=0$ can be written as\[
  \left(\frac{q_1}{q_0}\right)^2%
  +\left(\frac{q_2}{q_0}\right)^2%
  +\left(\frac{q_3}{q_0}\right)^2%
  +\left(\frac{q_4}{q_0}\right)^2%
  = -1.
\]
Now, given $a, b, c, d\in K$ with $a^2+b^2+c^2+d^2=-1$, it is tempting
to solve the equations $q_1=aq_0$, $q_2=bq_0$, $q_3=cq_0$, $q_4=dq_0$,
for $r,s,t,u$ in terms of $a,b,c,d$.

Indeed, this is possible, with
\begin{align*}
  r &= \frac{1}{3} c\\
  s &= \frac{a c^{2} d - 2 b c d}{2 a b^{2} c - b^{3} + b c^{2}}\\
  t &= \frac{-3 b d}{2 a b c - b^{2} + c^{2}}\\
  u &= \frac{-3 c d}{2 a b c - b^{2} + c^{2}}.
\end{align*}
Now with this choice of $r,s,t,u$, we first express $a_{14}, a_{17},
a_{21}, a_{22}$ in terms of $a,b,c,d$, and then solve for the
remaining $a_i$'s, obtaining Theorem \ref{theo}.
\section{Appendix}\label{app}
For completeness, we explain how a nowhere vanishing vector field on
the algebraic unit sphere produces a unimodular matrix as in Question
\ref{que}.

The following is a more elementary and down-to-earth version of
\cite[Remark 1.7]{AnaLev:combing}.

We use the notation from the introduction. Let $v(x, y, z)$ be a
nowhere vanishing vector field on $S(\bar K)$. Let
$V(X, Y, Z)\in K[X, Y, Z]^3$ be a preimage of $v(x, y, z)$.

First we see that $(x_0, y_0, z_0)$ and $V(x_0, y_0, z_0)$ are
linearly independent for all $(x_0,y_0,z_0)\in S(\bar K)$: Suppose
otherwise. Then $V(x_0,y_0,z_0)=\mu(x_0,y_0,z_0)$ for some
$\mu\in\bar K$. Then
$\mu=\mu
(x_0,y_0,z_0)\cdot(x_0,y_0,z_0)=V(x_0,y_0,z_0)\cdot(x_0,y_0,z_0)=0$, a
contradiction.

Write $V=(m_1, m_2, m_3)$. Then the three $2\times2$ minors of %
\[
  \begin{pmatrix}x_0 & y_0 & z_0\\m_1(x_0,y_0,z_0) & m_2(x_0,y_0,z_0)
    &m_3(x_0,y_0,z_0)\end{pmatrix}
\]
don't vanish simultaneously for any point $(x_0, y_0, z_0)\in S$. This
is the same as saying that the $2\times2$ minors of %
$\begin{pmatrix}X & Y & Z\\m_1 & m_2 &m_3\end{pmatrix}$ and
$X^2+Y^2+Z^2-1$ don't have a simultaneous solution in $\bar K^3$. Thus
by Hilbert's Nullstellensatz, there are
$m_4, m_5, m_6, m_7\in K[X, Y, Z]$ such that
\[
  m_4\cdot\smalldet{Y}{Z}{m_2}{m_3}%
  - m_5\cdot\smalldet{X}{Z}{m_1}{m_3}%
  + m_6\cdot\smalldet{X}{Y}{m_1}{m_2}%
  = m_7\cdot(X^2+Y^2+Z^2-1) + 1.
\]
As the left hand side is the determinant of
$M=\begin{pmatrix}X & Y & Z\\m_1 & m_2 & m_3\\m_4 & m_5 & m_6
\end{pmatrix}$, we obtain $\det(M)\equiv1\pmod{X^2+Y^2+Z^2-1}$, hence
$\det(M(x, y, z))=1$.
\printbibliography
\end{document}